%% file: LeastSquares.tex
\documentclass[10pt]{article}
\usepackage{comment}
\usepackage{fancybox}
\usepackage{color}
\definecolor{light-gray}{gray}{0.9}

\date{}
\include{preamble}

\title{`Basic' Generalization Error Bounds for  Least Squares Regression with Well-specified Models}
\author{Karthik Duraisamy, \\ University of Michigan, Ann Arbor}

\input{macros}

\begin{document}
\maketitle

\begin{abstract}
	This note examines the behavior of generalization capabilities - as defined by out-of-sample mean squared error (MSE)  - of Linear Gaussian (with a fixed design matrix) and Linear Least Squares regression. Particularly, we consider a well-specified model setting, i.e. we assume that there exists a `true'  combination of model parameters within the chosen model form. While the statistical properties of Least Squares regression have been extensively studied over the past few decades - particularly with {\bf less restrictive problem statements} compared to the present work -  this note targets bounds that are {\bf non-asymptotic and more quantitative} compared to the literature. Further, the analytical formulae for distributions and bounds (on the MSE) are  directly compared to numerical experiments.  Derivations are presented in a self-contained  and pedagogical manner, in a way that a reader with a basic knowledge of probability and statistics can follow.
\end{abstract}
\vspace{1cm}

\section{Introduction}

We consider a linear model $p(y|x,\theta) = \mathcal{N}(y;\phi(x)^T \theta, \sigma^2)$, where  $y \in \reals$ is the response function, $x \in \reals^p$ is the input,  $\theta \in \reals^m$ represents model parameters, and $\phi : \reals \rightarrow \reals^m$ is the feature vector. Given a dataset $D =\{(x_1,y_1) ... (x_n,y_n)\}$, where the inputs are drawn in an i.i.d. fashion from a known $p(x)$, define $\Phi(x) \in \reals^{n \times m}$, where each row of $\Phi(x)$ is $\phi(x_i)^T$. Further, we assume that the data is generated by a true set of parameters $\theta^\ast \in \reals^m$ , and that the noise in the data is independent of the inputs $x_i$.  Finally, we assume that $\phi(x)$ follows the standard normal distribution.

 Define a squared objective $\ell (\theta,x,y) \triangleq \sum_{i=1}^n (y_i - \phi(x_i)^T \theta)^2$. Assuming $n>m$, and that $\Phi$ has full-rank,  minimization of the objective function over the dataset yields the least squares estimate $\theta^{LS} = (\Phi^T(x) \Phi(x))^{-1} \Phi(x)^T y$. Now, given a new input $\xh \in \reals^p$, we would like to predict the response and assess the loss  ${\ell}(\theta,\hat{x},\hat{y}) = (\hat{y}-\phi(\hat{x})^T \theta)^2,$ where $\yh$ is a realization from $\mathcal{N}(\phi^T(\xh)\theta^{\ast},\sigma^2)$. 

Note that the above setting corresponds to random design regression with a well-specified linear model.  Breiman and Freedman~\cite{breiman1983many}, showed that~\footnote{An additional nuance of their contribution is pointed out in section~\ref{sec:LSQ}}
\begin{align}
\mathbb{E}[\ell (\theta^{LS},\xh,\yh)] = \sigma^2  \frac{n-1}{n-1-m}.
\end{align}
We are interested in studying the statistical properties of ${\ell}(\theta^{LS},\hat{x},\hat{y}) $  for out-of-sample predictions, particularly in obtaining bounds on $\ell(\cdot,\cdot,\cdot)$.

Indeed, statistical properties and generalization of least squares regression have long been studied, including in less restrictive settings (e.g. mis-specified models, non-linear features, non-Gaussian distributions, etc.). However, to the knowledge of the author~\footnote{The author is willing to expand his understanding of the literature, and would be glad to be proved wrong}, much of the literature - though rigorous, well-crafted and require {
\em fewer assumptions} than the above - presents bounds that are {\em  more restrictive} in a number of ways. The following are examples:

$\bullet$ A number of works present bounds for high dimensional and /or asymptotic situations  (e.g. ~\cite{hastie2019surprises}) or require bounded covariates.

$\bullet$ Gyorfi et al.~\cite{gyorfi2002distribution} and Catoni's~\cite{catoni2004statistical} error bounds are well-crafted, but contain an arbitrary constant. 

$\bullet$ Many publications contain terms such as $O(\cdot)$, thus rendering an asymptotic error estimate.

$\bullet$ Audibert \& Catoni~\cite{audibert2010linear} require 
$n >> m \log m$ and Hsu et al.~\cite{hsu2012random} require $n >> m$.

$\bullet$ Classical PAC-Bayesian bounds~\cite{mcallester1999pac} require bounded loss functions or additional parameters beyond the data~\cite{germain2016pac}.

It is also intriguing that in much of the literature, bounds are often not verified using numerical experiments (indeed, there are clear exceptions, for instance~\cite{hastie2019surprises,mei2019generalization}). The author is careful to re-emphasize that many of the above works focus on a more general - and thus more practically relevant - setting than the present one, and that a bound can be useful even under the above conditions.

The main contribution of the present work is the following: Assuming the above conditions and $n>m-3$, with a probability of $1-\delta$, where $\delta \in (0,1)$, we can guarantee that

\begin{equation*}
\ell(\xh,\yh,\theta^{LS}) \leq  \frac{\sigma^2}{n-m-1} \left[m+\frac{1}{\sqrt \delta}\sqrt{\frac{m^3 - m^2 (n-3) - 3 m (n-3) (n-1) + 
		3 (n-3) (n-1)^2}{n-m-3} } \right].
\end{equation*}
It is noted that this result does {\bf not} require asymptotics.  

The plan is to proceed with the derivation in a pedagogical manner, in such a way that a reader with the most basic knowledge of probability can follow. Given the above goal in combination with the restrictions placed by the assumptions in the problem statement, prompted the use of the words ``Basic Bounds". Aligned with the pedagogical spirit, we first start with the linear Gaussian setting (i.e. The design matrix is a constant) and derive the {\bf distributions} of the training and testing (generalization) error. Then we develop {\bf bounds} for the generalization error in the linear Gaussian and  linear least squares regression settings.

\section{Linear Gaussian with Least Squares}

\begin{theorem}[Well-specified Linear Gaussian case]\label{th:lingauss}
Consider $p(y|\theta) = \mathcal{N}(A \theta, \sigma^2 I)$ with $y\in \reals^n$, $\theta \in \reals^m$, and $A \in \reals^{n\times m}$ being a full column rank matrix . Define the empirical risk as a log-likelihood-type loss
\begin{eqnarray}
{\ell}(y,\theta) \triangleq     \frac{1}{2 \sigma^2}   (y-A\theta)^T (y-A\theta).
\end{eqnarray}
Further, assume that the data is generated by a true process $\Theta^\ast \sim \mathcal{N}(\theta^\ast,0).$ 

Consider a training sample $y \triangleq A\theta^\ast + \sigma^2 z$, and a testing sample $y \triangleq A \theta^\ast + \sigma z_t$, where $Z,Z_t \sim \mathcal{N}(0,I_n)$. Define the least squares estimate $\theta^{LS} \triangleq  (A^T A)^{-1} A^T y$.  Then the  training,  testing, and true risks are given by 
	\begin{align}
	\mathcal{L}(Y,\Theta^{LS})  &\sim \frac{1}{2} \left[ \chi^2(n-m) \right]; \ \ \ \ \mathbb{E}[\mathcal{L}(Y,\Theta^{LS}) ] = \frac{1}{2} (n-m) \\
		\mathcal{L}(Y_t,\Theta^{LS})  &\sim \frac{1}{2} \left[  \chi^2(n-m) + 2 \chi^2(m) \right] ; \ \ \mathbb{E}[\mathcal{L}(Y_t,\Theta^{LS}) ] = \frac{1}{2} (n+m)\\
			\mathcal{L}(Y,\Theta^\ast)  &\sim \frac{1}{2} \left[ \chi^2(n) \right]; \ \ \ \ \mathbb{E}[\mathcal{L}(Y,\Theta^\ast) ] = \frac{1}{2} n 
	\end{align}
Additionally, given $0 \leq \delta \leq 1$, and a random testing sample $y_t$, with at least a probability of $1-\delta$, one can guarantee that the generalization error is bounded by

	\begin{align}
	\norm{y_t-A\theta^{LS}}^2 \leq  n+m + \sqrt{\frac{6m+2n}{\delta}}.
	\end{align}

\end{theorem}

{\bf Proof:}



Let's consider the training risk first:
\begin{align*}
    (y-A\theta)^T (y-A\theta) &=   (A\theta^\ast + \sigma z - AA^+ (A \theta^\ast +\sigma z ) )^T (A\theta^\ast + \sigma z - AA^+ (A \theta^\ast +\sigma z )) \\
&=  (\sigma z - \sigma P z  )^T (\sigma z - \sigma P z) \\
&=  \sigma^2 z^T P z  + \sigma^2  z^T z -  2 \sigma^2 z^T P z, \\
\end{align*}
where $P \triangleq A A^+ = A^{+T} A^T,$ is the projection matrix and is of rank $r$. Since A has full column rank, $r=m$.
Therefore  
\begin{align*}
{\ell}(y,\theta^{LS}) &= \frac{1}{2 \sigma^2} \left[ \sigma^2 z^T z-   \sigma^2  z^T P z \right] =  \frac{1}{2} \left[ z^T (I-P)z \right].
\end{align*}

This is a random variable with the following Chi-squared distribution~\footnote{See the testing risk derivation for more insight into how the structure of P is leveraged}
\begin{align}
\mathcal{ L}(\Theta^{LS})  &\sim \frac{1}{2} \left[\chi^2(n-r) \right] \ \ = \Gamma \left(\frac{n-m}{2},1 \right),
\end{align}
where $\Gamma$ is the Gamma distribution.
The expected empirical risk is 

\begin{align}
\mathbb{E} [\mathcal{ L}(\Theta^{LS})]  &=  \frac{1}{2}  ( n -  m ).
\end{align}

Similarly, the true risk is 
\begin{align}
\mathcal{ L}(\Theta^\ast) & \sim \frac{1}{2 \sigma^2}\left[\sigma^2 \chi^2(n)\right]\ \ =   \Gamma \left(\frac{n}{2},1 \right).
\end{align}

The expected true risk is therefore
\begin{align}
\mathbb{E} [\mathcal{ L}(\Theta^\ast)]  &=  \frac{n}{2}.   
\end{align}

\noindent {\bf Testing risk}

Let's examine the testing risk with $y_t = A \theta^\ast + \sigma z_t$, and $\theta = \theta^\ast + \sigma A^+  z$   where $Z_t, Z\sim \mathcal{N}(0,I_n)$ represent testing and training realizations. Then

\begin{align*}
(y_t-A\theta)^T (y_t-A\theta) &=   (\sigma z_t -  \sigma P z)^T (\sigma z_t- \sigma P z) \\
&=  \sigma^2 (z^T P z  +  z_t^T z_t -  2  z^T P z_t).
\end{align*}
Write $P = S \Lambda S^T$, where $S$ is an orthonormal matrix. Also define $Q \triangleq S^T Z$ and $Q_t \triangleq S^T Z_t$.
Define $ g \triangleq z^T P z  +  z_t^T z_t -  2  z^T P z_t$.
\begin{align*}
g &= q^T \Lambda q + q_t^T q_t - 2 q^T \Lambda q_t \\
 &=  \sum_{j=1}^n ( \lambda_j q_j^2 + q_{tj}^2 - 2 \lambda_j q_j q_{tj} ) \\
  &=  \sum_{j=1}^r ( q_j^2 + q_{tj}^2 - 2 q_j q_{tj} ) + \sum_{j=r+1}^n q_{tj}^2  \\
    &=  \sum_{j=1}^r ( q_j - q_{tj} )^2 + \sum_{j=r+1}^n q_{tj}^2.
\end{align*}
The fact that $P$ has exactly $r$ unity eigenvalues and $n-r$ zero eigenvalues has been utilized above. Recognizing that $Q_j-Q_{tj} \sim \mathcal{N}(0,2)$, it is easy to see that
  $  G \sim 2 \chi^2(r) + \chi^2(n-r)$.

The testing risk is therefore

$$ \mathcal{L}(\theta^{LS})  =\frac{1}{2 \sigma^2} \left[ \sigma^2  (2 \chi^2(r) + \chi^2(n-r)) \right].$$ Again, since $A$ has full column rank, $r=m$, and thus

\begin{empheq}[box=\tcbhighmath]{equation*}
	\mathcal{L}(\theta^{LS})  =\frac{1}{2} \left[  2 \chi^2(m) + \chi^2(n-m) \right].
\end{empheq}

Note: One has to be careful in handling the above expression because of the coefficients in front of the $\chi^2(\cdot)$ distributions are different.  Nevertheless, the expected test risk can be evaluated easily:

\begin{align}
\mathbb{E} [\mathcal{L}(\theta^{LS})]  &=  \frac{n+m}{2}.
\end{align}

\subsection{Numerical Verification}

\begin{figure}
	\centering
	\includegraphics[width=0.8\textwidth]{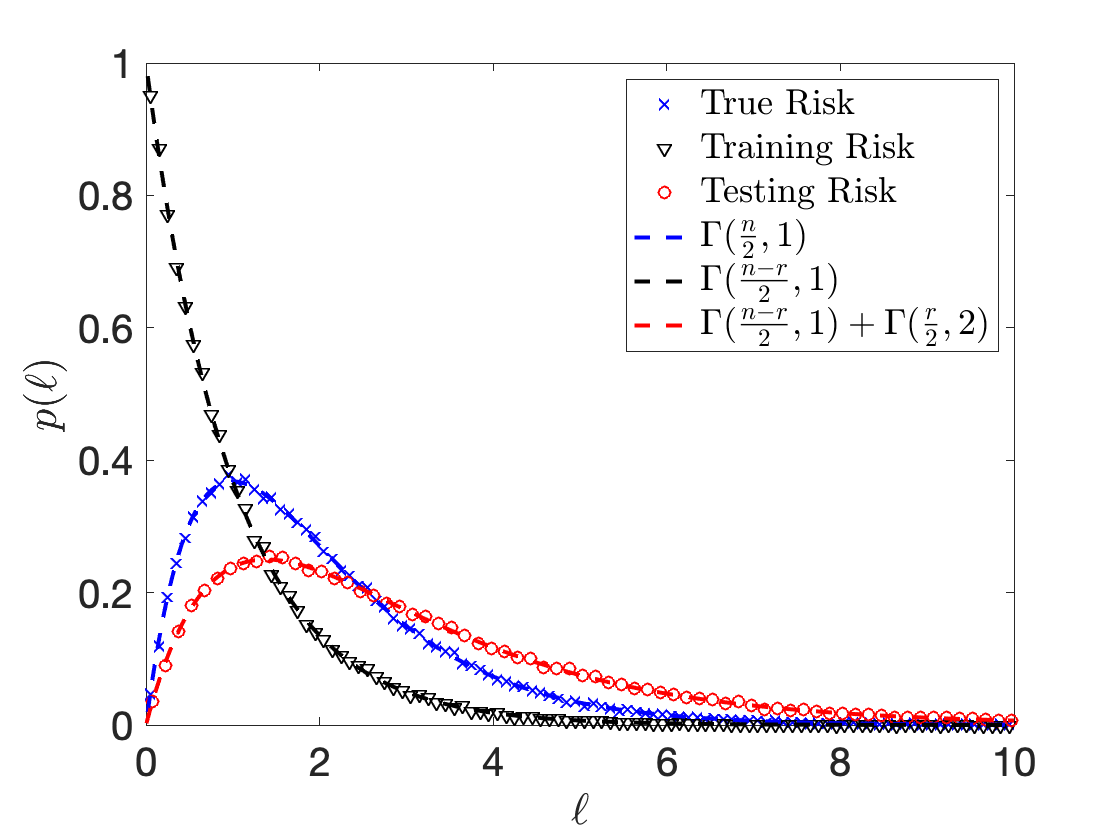}
	\caption{PDF of risks  determined by Sampling (Symbols) compared to predictions (lines).}
	\label{fig:Risk}
\end{figure}

To verify the above equations, a sample problem is designed with 
$$ \sigma = 0.1 \ \ ; \ \ \theta^\ast = \left[ 
\begin{array}{c }
0.3 \\
-2 \\
\end{array}
\right] \  \ ; \  \
A= \left[
\begin{array}{c c}
1 & 0.6 \\
3.2 & -2 \\
4 & 1 \\
3.1 & -1 \\
\end{array}
\right].
$$
and the results are shown in Figure~\ref{fig:Risk}. A total of 100,000 samples were used to evaluate the 3 different risk elements which are compared to the analytical forms.

\begin{table}[htbp]
	{   
		\caption{Risk Decomposition for Least Squares}
		\label{tab:density}
		\begin{center}
			\begin{tabular}{|c|c|c|l|l|} \hline
			\bf Term &	\bf Expression & \bf Distribution  & \bf Expectation  & \bf Variance \\ \hline \hline
			Training Risk		&	$ \norm{Y-A \Theta^{LS}}^2 $ & $\sigma^2 \chi^2(n-m)$   &   $  \sigma^2( n-m) $ & $\sigma^4 2 (n-m) $ \\ \hline
			True Risk &	$ \norm{Y-A \Theta^\ast}^2 $ & $\sigma^2 \chi^2(n)$   &   $  \sigma^2( n ) $  & $\sigma^4 2 n $ \\ \hline
			Testing Risk &	$ \norm{Y_t-A \Theta^{LS}}^2 $ & $\sigma^2 (\chi^2(n-m) + 2 \chi^2(m))$   &   $  \sigma^2( n+m) $ & $\sigma^4 (6m + 2n) $ \\ \hline
			\end{tabular}
		\end{center}
	}
\end{table}

\subsection{Bounds}

\begin{figure}
	\centering
	\includegraphics[width=0.49\textwidth]{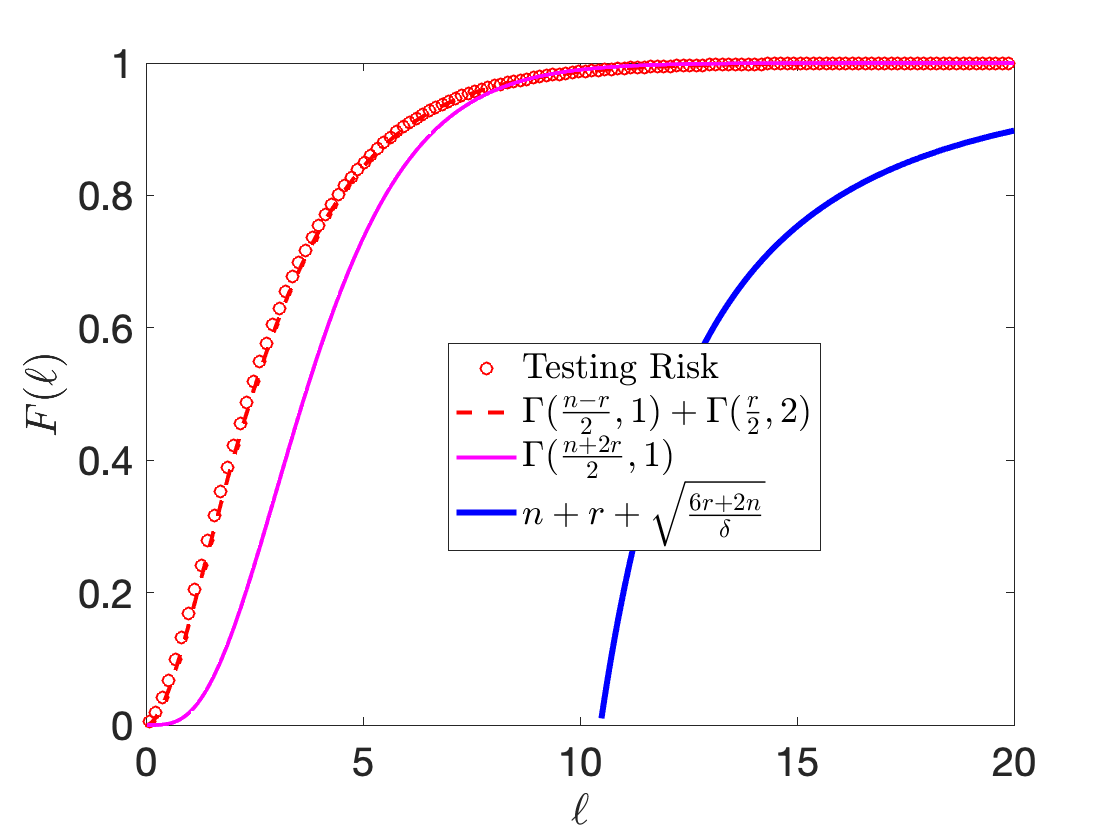}
	\includegraphics[width=0.49\textwidth]{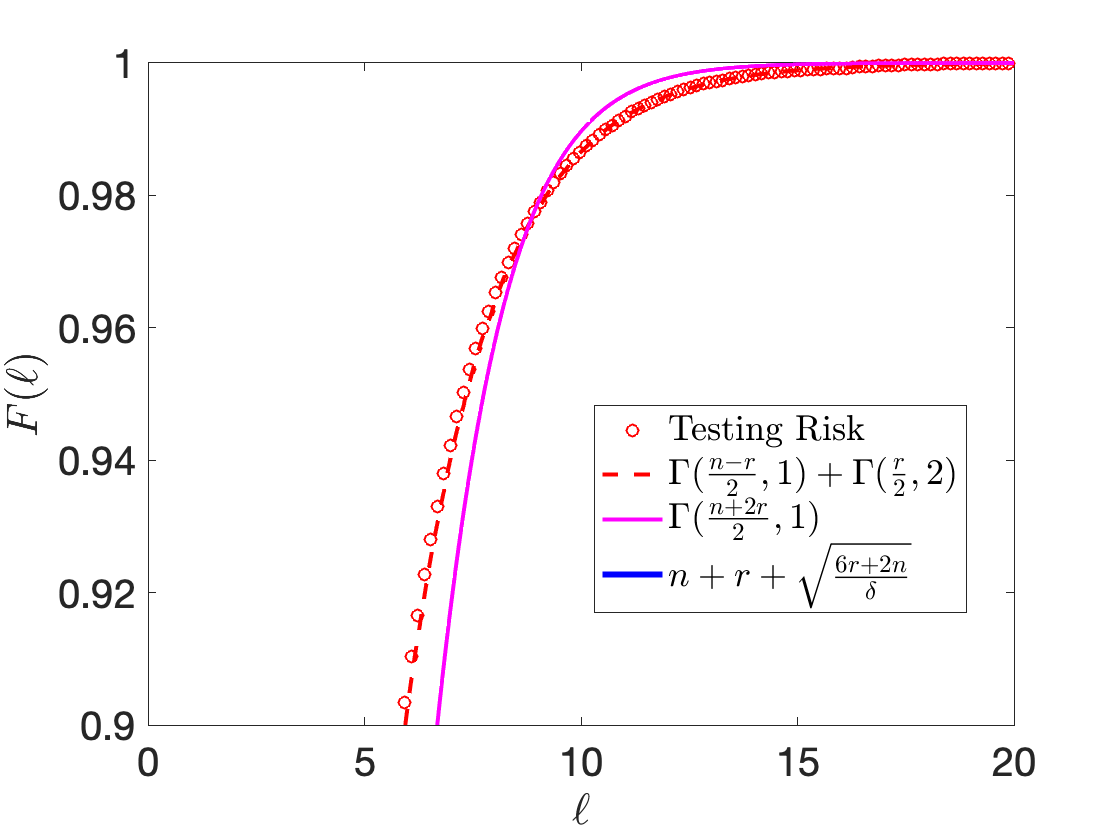}
	\caption{CDF of testing risk  determined by Sampling (Symbols) compared to estimations (lines).}
	\label{fig:Risk1}
\end{figure}

Let's try to bound the test risk.
It is perhaps intuitive to  assume that $2 \chi^2(m) + \chi^2(n-m)$ can be approximated as $\chi^2(n+m)$. However, that is highly inaccurate, and neither does the more conservative approximation of $\chi^2(n+2m)$ serve as a good bound as seen in Figure~\ref{fig:Risk1}.

We can, however, use Chebyshev's inequality, which states that for any random variable $X$, and $\epsilon \in \reals^+$,

$$ P( |x-\mathbb{E}[X]| \geq \epsilon)  \leq  \frac{Var(X)}{\epsilon^2}. $$

Therefore

\begin{align*}
P( |g-(n+m)| \geq \epsilon) &  \leq  \frac{6m+2n}{\epsilon^2} \\
P( |g-(n+m)| \leq \epsilon)   &\geq 1- \frac{6m+2n}{\epsilon^2} \\
\end{align*}
Define $\delta \triangleq \frac{6m+2n}{\epsilon^2} $. Therefore, $\epsilon = \sqrt{\frac{6m+2n}{\delta}}$.

Therefore, with at least a probability of $1-\delta$, we can guarantee that 
$$
 |g-(n+m)| \leq  \sqrt{\frac{6m+2n}{\delta}}
$$

Therefore with at least a probability of $1-\delta$, we can guarantee that the generalization error is bounded by

\begin{empheq}[box=\tcbhighmath]{equation*}
\norm{y-A\theta^{LS}}^2 \leq  n+m + \sqrt{\frac{6m+2n}{\delta}}.
\end{empheq}

This expression is verified in Figure~\ref{fig:Risk1}.


Better bounds may be obtained with higher moments and concentration inequalities.


\section{Linear Least Squares Regression}
\label{sec:LSQ}

\begin{theorem}[Well-specified Linear Least Squares Regression]\label{th:linleast}
Consider  $p(y|x,\theta^\ast) = \mathcal{N}(\phi(x)^T \theta^\ast, \sigma^2)$, where  $y \in \reals$, $x \in \reals^p$, and $\theta^\ast \in \reals^m$. $\phi : \reals \rightarrow \reals^m$ represents the feature vector. We will also assume that $\Phi(x)$ has full column rank.
Given a dataset $D =\{(x_1,y_1) ... (x_n,y_n)\}$, define $\Phi(x) \in \reals^{n \times m}$, where each row of $\Phi(x)$ is $\phi(x_i)^T$.
Let's define a MSE loss ${\ell}(\theta,\hat{x},\hat{y}) \triangleq (\hat{y}-\phi(\hat{x})^T \theta)^2$, and the OLS estimate  $\theta^{LS}(y) \triangleq (\Phi(x)^T \Phi(x))^{-1} \Phi(x)^T y = \Phi^+(x) y$. 

Then, for a random $\hat{x}$, assuming that the features are distributed in a Gaussian distribution (i.e.  $\phi(x) \sim \mathcal{N}(0,\Sigma_\phi))$, with at least a probability of $1-\delta$ the generalization error is bounded by 
\begin{equation*}
\ell(\theta^{LS},\hat{x},\hat{y}) \leq  \frac{\sigma^2}{n-m-1} \left[m+\frac{1}{\sqrt \delta}\sqrt{\frac{m^3 - m^2 (n-3) - 3 m (n-3) (n-1) + 
		3 (n-3) (n-1)^2}{n-m-3} } \right],
\end{equation*}
with $n>m-3$.

\end{theorem}

{\bf Proof:} 
 In the below, we will use  $\theta$ instead of $\theta^{LS}$ to improve readability.

\begin{align*}
{\ell}(\theta,\hat{x},\hat{y}) &= \norm{\phi(\hat{x})^T \theta^\ast + \sigma \zh - \phi(\hat{x})^T \theta }^2\\
&= \norm{\phi(\hat{x})^T (\theta^\ast - \theta)+ \sigma \zh }^2 \\
&= \norm{\phi(\hat{x})^T (\theta^\ast - \Phi^+(x) y)+ \sigma \zh }^2 \\
&= \norm{\phi(\hat{x})^T (\theta^\ast - \Phi^+(x) (\Phi(x) \theta^\ast + \sigma z))+ \sigma \zh }^2 \\
&= \sigma^2 \norm{\zh -\phi(\hat{x})^T \Phi^+(x) z  }^2 \\
\mathbb{E}_{\zh}[{\ell}(\theta,\hat{x},\hat{y})]&=\sigma^2+ \sigma^2 \norm{\phi(\hat{x})^T \Phi^+(x) z}^2  \\
\mathbb{E}_{z} \mathbb{E}_{\zh} [{\ell}(\theta,\hat{x},\hat{y})]&=\sigma^2+ \sigma^2 \tr[\phi(\hat{x})^T \Phi(x)^+ \Phi^{+T}(x) \phi(x)] \\
&= \sigma^2+\sigma^2 \tr[\phi(\hat{x})^T (\Phi^{T}(x) \Phi(x))^{-1} \phi(x)] \\
\mathbb{E}_{\xh}\mathbb{E}_{z} \mathbb{E}_{\zh} [{\ell}(\theta,\hat{x},\hat{y})]&=\sigma^2+ \sigma^2 \tr[\Sigma_\phi  (\Phi(x)^T \Phi(x))^{-1}],
\end{align*}

where $\Sigma_\phi = \mathbb{E}[\phi(\xh) \phi(\xh^T)]$. Note that the above quantity is still a random variable because of the randomness in $x$.

\noindent {\bf Gaussian assumption for features:} Assume  $\phi(x) \sim \mathcal{N}(0,\Sigma_\phi)$. Writing the Cholesky decomposition $\Sigma_\phi = L L^T$, and thus $\phi(x_i) = L q_i$, where $q_i \sim \mathcal{N}(0,I_m)$. Note: we are temporarily dropping small / big symbols for realizations and random variables to stop the proliferation of symbols.
$$ \Phi(x) = \left[ 
\begin{array}{c }
\phi(x_1)^T \\
..\\
..\\
\phi(x_n)^T \\
\end{array}
\right]  = 
\left[ 
\begin{array}{c }
q_1^T  L^T \\
..\\
..\\
q_n^T  L^T \\
\end{array}
\right] 
= 
\left[ 
\begin{array}{c }
q_1^T   \\
..\\
..\\
q_n^T  \\
\end{array}
\right]  L^T \triangleq Q^T L^T
$$

Therefore $\Phi(x)^T \Phi(x) = L Q Q^T L^T$, where $Q \in \reals^{m \times n}.$ Thus 

\begin{align*}
\mathbb{E}_{\xh}\mathbb{E}_{z} \mathbb{E}_{\zh} [{\ell}(\theta,\hat{x},\hat{y})]&= \sigma^2+\sigma^2 \tr[L L^T (L Q Q^T L^T)^{-1}]\\
&= \sigma^2+\sigma^2 \tr[L L^T L^{-T}(Q Q^T)^{-1}L^{-1}]\\
&=\sigma^2+ \sigma^2 \tr[(Q Q^T)^{-1}]\\
\mathbb{E}_x\mathbb{E}_{\xh}\mathbb{E}_{z} \mathbb{E}_{\zh} [{\ell}(\theta,\hat{x},\hat{y})] &= \sigma^2+\sigma^2 \tr[\mathbb{E}((Q Q^T)^{-1})].
\end{align*}

We know that $QQ^T \sim \mathcal{W}_m(I_m,n)$, and $(QQ^T)^{-1} \sim \mathcal{W}^{-1}_m(I_m,n)$, where $\mathcal{W}_m$ and $\mathcal{W}_m^{-1}$ represent the Wishart and Inverse Wishart~\cite{gelman2013bayesian} distributions, respectively. Therefore $\mathbb{E}((Q Q^T)^{-1}) = \frac{I_m}{n-m-1}$ as long as $n > m+1$.

Therefore 
\begin{align}
\mathbb{E}[{\ell}(\theta,\hat{x},\hat{y})] \triangleq \mu_\ell = \sigma^2+ \sigma^2 \frac{m}{n-m-1} .
\end{align}
This expression can also be found in Breiman \& Freedman~\cite{breiman1983many}. We are interested in bounds of ${\ell}(\theta,\hat{x},\hat{y})$, and  note that Theorem 1.3 of Ref.~\cite{breiman1983many}, instead provides the distribution of $\mathbb{E}[\ell(\theta,\hat{x},\hat{y})|x,y]$.
Note that the true expected risk is $\mathbb{E}[{\ell}(\theta^\ast,\hat{x},\hat{y})] = \sigma^2.$

%
%

\subsection{Variance}
Let's try to bound the above distribution by computing the variance of the squared error
\begin{align*}
{\ell}(\theta,\hat{x},\hat{y})^2 
&= \sigma^4(\zh^T-z^T\Phi^{+T} \phi )^4 .
\end{align*}

We will use the following identity~\cite{matrixref}:
If $X \sim \mathcal{N}(0,I)$, $A$ is a matrix and $a$ is a vector, then
$$\mathbb{E}_X[(Ax + a)^T(Ax + a) (Ax + a)^T(Ax + a)] = 2\tr(AA^TAA^T) + 4a^TAA^Ta + (tr(AA^T) + a^T a)^2.$$

Thus, 
\begin{align*}
\mathbb{E}_{\hat z} \ell^2 &= \sigma^4(2 + 4 (z^T\Phi^{+T} \phi)^2 + (1+(z^T\Phi^{+T} \phi)^2 )^2)\\
&= \sigma^4(3 + 6 (\phi^T  \Phi^{+} z)^2 + (\phi^T \Phi^{+} z)^4)  \\
\mathbb{E}_{z} \mathbb{E}_{\hat z} \ell^2 
 &= \sigma^4(3 + 6 \tr(\Phi^{+T}\phi \phi^T \Phi^+) + 2 \tr(\phi^T \Phi^{+}  \Phi^{+T}\phi \phi^T \Phi^{+}  \Phi^{+T}\phi)+ \tr(\phi^T \Phi^{+}  \Phi^{+T}\phi)^2 ) \\
 &= \sigma^4(3 + 6 \tr(\Phi^{+T}\phi \phi^T \Phi^+) + 3 \phi^T \Phi^{+}  \Phi^{+T}\phi \phi^T \Phi^{+}  \Phi^{+T}\phi ) \\ 
\mathbb{E}_{\xh} \mathbb{E}_{z} \mathbb{E}_{\hat z} \ell^2&= \sigma^4(3 + 6 \tr(\Phi^{+T} \Phi^+) + 6 \tr(\Phi^{+T}  \Phi^{+} \Phi^{+T}  \Phi^{+}) + 3 \tr(\Phi^{+T}  \Phi^{+})^2) \\
&= \sigma^4(3 + 6 \tr( (\Phi^{T} \Phi)^{-1}) + 6 \tr((\Phi^{T} \Phi)^{-2}) + 3 \tr((\Phi^{T} \Phi)^{-1})\tr((\Phi^{T} \Phi)^{-1}))\\
&\sim \sigma^4(3 + 6 \tr( W_m^{-1}(I,n)) + 6 \tr(W_m^{-1}(I,n)W_m^{-1}(I,n)) + 3 \tr(W_m^{-1}(I,n))\tr(W_m^{-1}(I,n))),
\end{align*}
where $W_m^{-1}$ represents the Inverse Wishart distribution. The expectation of the last expression can be compactly reduced using standard Wishart distribution identities, except the term involving the trace of the product of inverse Wishart matrices, which we obtain from Pielaszkiewicz \& Holgersson~\cite{pielaszkiewicz2019mixtures} (page 8). With this identity,
\begin{align*}
\mathbb{E}_x \mathbb{E}_{\xh} \mathbb{E}_{z} \mathbb{E}_{\hat z} \ell^2&= \sigma^4\left(3 + 6 \frac{m}{n-m-1}  + 6 \frac{(n-1)m}{(n-m-3)(n-m-1)(n-m)}+ 3  \frac{m(m(n-m-2)+2)}{(n-m-3)(n-m-1)(n-m)}\right) \\
&=\sigma^4 \frac{3 (n-1) (n-3)}{(n-m-1) (n-m-3)},
\end{align*}
as long as $n>m+3$.
Therefore, the variance is

$$\sigma^2_\ell = \mathbb{E}[\ell^2]- \mathbb{E}[\ell]^2 = \sigma^4 \frac{m^3 - m^2 (n-3) - 3 m (n-3) (n-1) + 
3 (n-3) (n-1)^2}{(n-m-1)^2 (n-m-3)}.
$$

\begin{figure}
	\centering
	\includegraphics[width=0.49\textwidth]{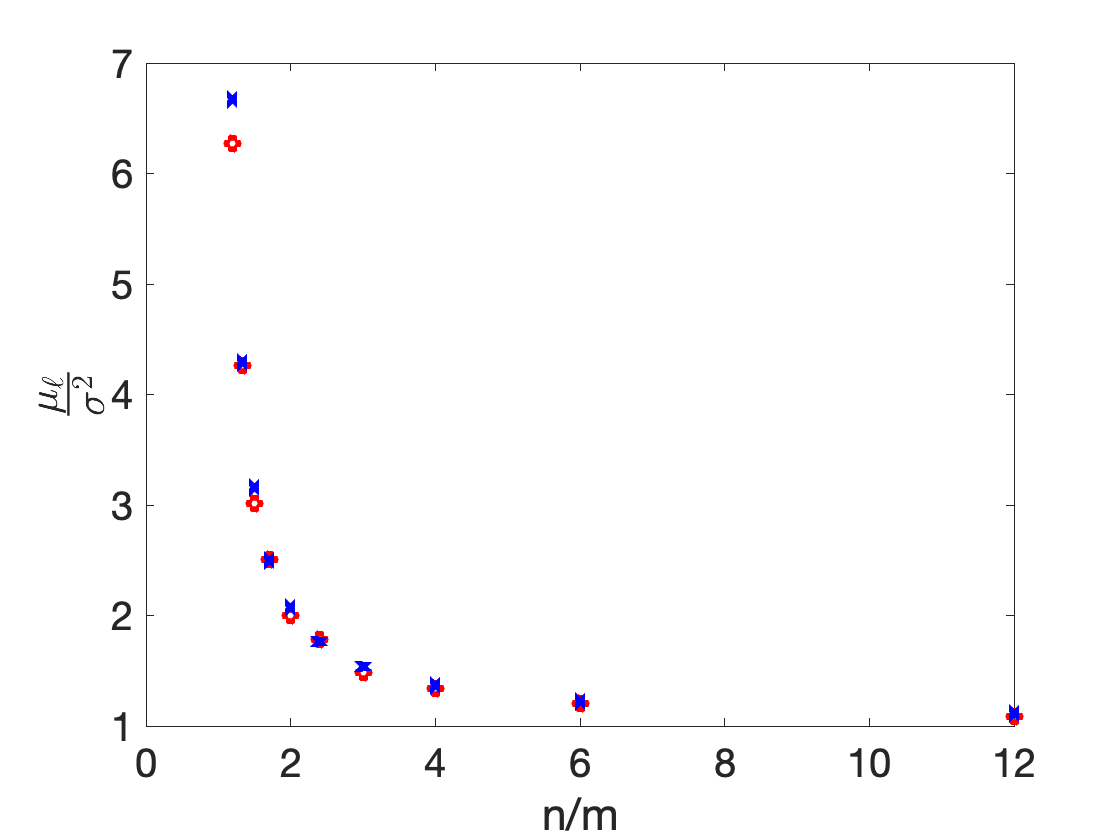}
	\includegraphics[width=0.49\textwidth]{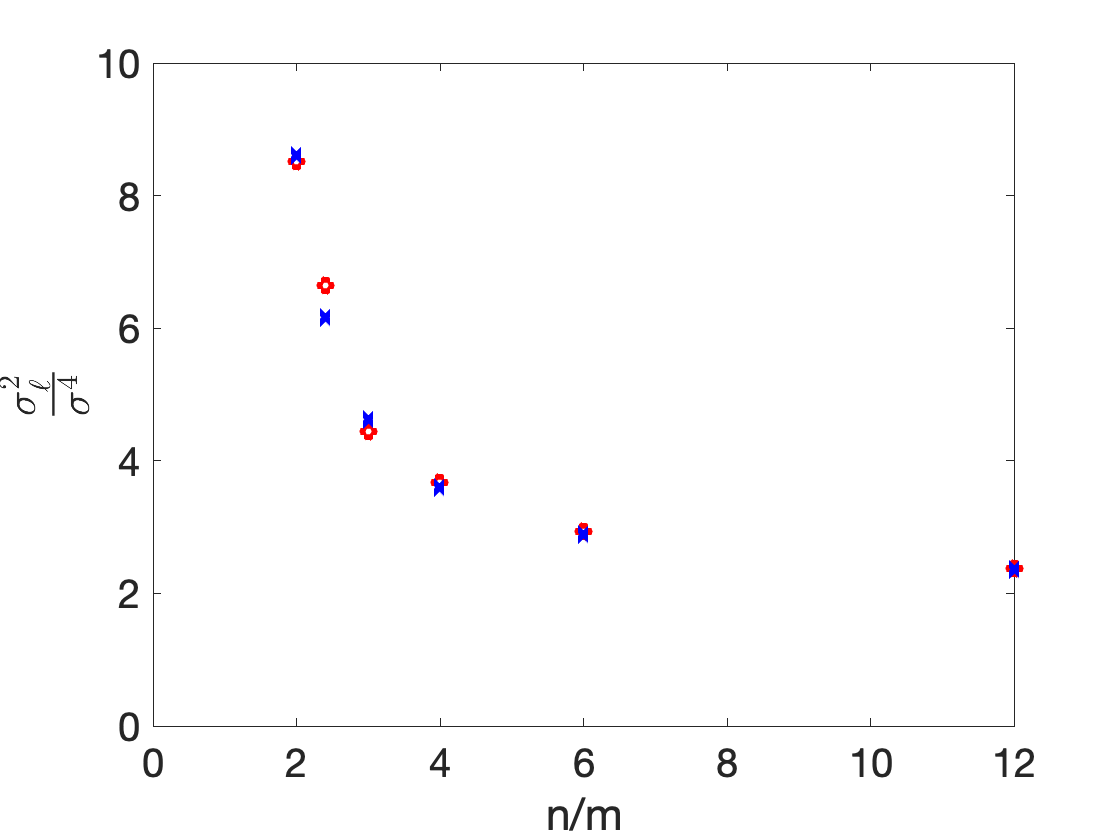}
	\caption{Analytical vs empirical mean and variance of MSE for OLS with $n=60$. Empirical means were computed using 100 different simulations, with the mean and variance evaluated using $100n$ samples (for each experiment. Thus, a total of 10000n samples)}
	\label{fig:MSEOLS}
\end{figure}

The analytical mean and variance are cmpared to numerically evaluated ones in Figure~\ref{fig:MSEOLS}. The numerical experiments use $X \sim \mathcal{N}(0,1)$, $\sigma^2 = 0.04$, $\phi(x)=x$, and $\theta^\ast = \frac{t}{\norm{t}}$, where $T \sim \mathcal{N}(0,I_m)$.

To obtain  bounds, we employ the Chebyshev inequality
\begin{align*}
P( |\ell-\mu_\ell| \geq \epsilon) &  \leq  \frac{\sigma_\ell^2}{\epsilon^2} \\
P( |\ell-\mu_\ell| \leq \epsilon)   &\geq 1- \frac{\sigma_\ell^2}{\epsilon^2}.
\end{align*}
Define $\delta \triangleq \frac{\sigma_\ell^2}{\epsilon^2} $. Therefore, $\epsilon = {\frac{\sigma_\ell}{\sqrt \delta}}$.

Therefore with at least a probability of $1-\delta$, we can guarantee that 
$$
|\ell-\mu_\ell| \leq \frac{\sigma_\ell}{\sqrt \delta}.
$$
The upper tail bound for the squared error  with a probability of $1-\delta$ is

\begin{empheq}[box=\tcbhighmath]{equation*}
\ell \leq  \frac{\sigma^2}{n-m-1} \left[m+\frac{1}{\sqrt \delta}\sqrt{\frac{m^3 - m^2 (n-3) - 3 m (n-3) (n-1) + 
		3 (n-3) (n-1)^2}{n-m-3} } \right].
\end{empheq}

\begin{figure}
	\centering
	\includegraphics[width=0.49\textwidth]{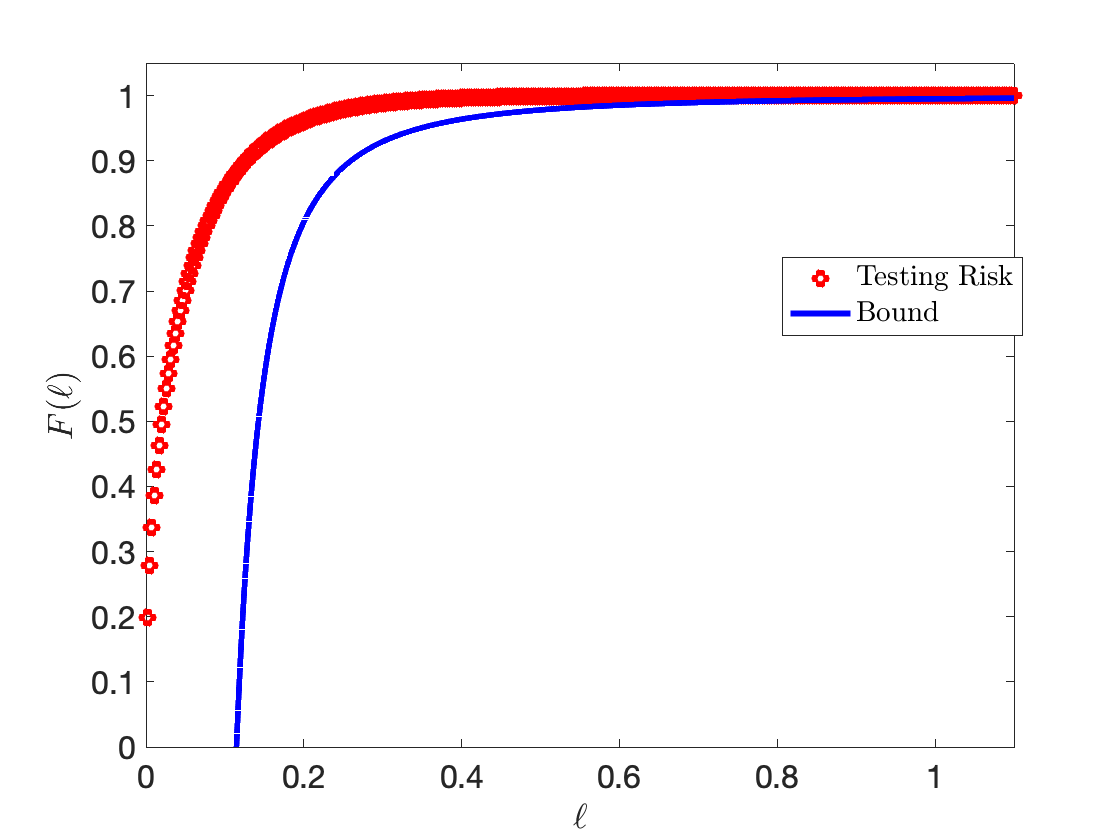}
	\includegraphics[width=0.49\textwidth]{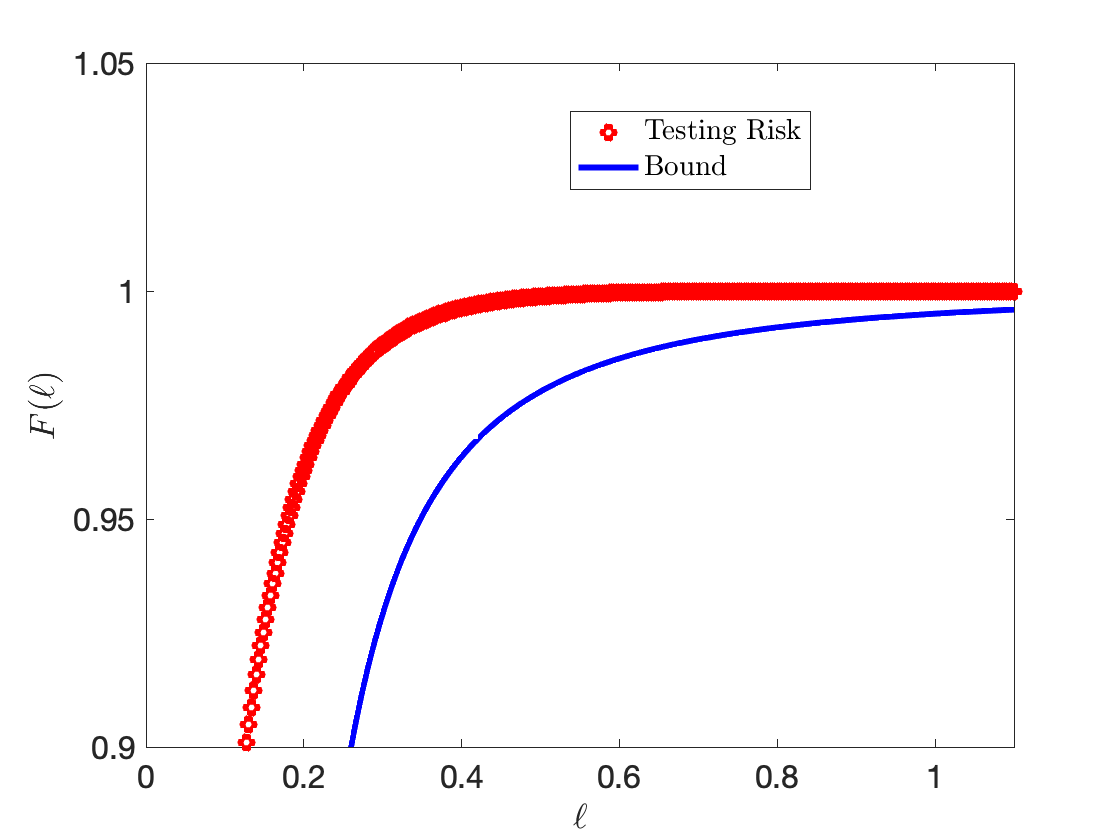}
	\caption{Empirical CDF and tail bounds of OLS with $n=60$, $m=10$.}	
	\label{fig:TailOLS}
\end{figure}

This bound is evaluated using the afore-mentioned numerical experiments in Figure~\ref{fig:TailOLS}. 

\subsection{Finite, but large n}

For $n>>1$ but still finite, we have

$$\sigma^2_\ell \approx \sigma^4 \frac{m^3 - m^2 n - 3 m n^2 + 
	3 n^3}{(n-m)^3 }.  $$

defining $\alpha \triangleq n/m$ and dividing by $m^3$

$$\sigma^2_\ell \approx \sigma^4 \frac{1 - \alpha - 3\alpha^2 + 
	3 \alpha^3}{(\alpha-1)^3 } = \sigma^4 \frac{ 3\alpha^2-1}{(\alpha-1)^2 }.$$

Thus for finite $n >>1$ we have 
$$
\mu_\ell \approx \sigma^2+ \frac{\sigma^2}{\alpha-1} = \frac{\alpha}{\alpha-1}\sigma^2 \ \ ; \ \  \sigma_\ell \approx \frac{\sigma^2}{\alpha -1} \sqrt{ 3\alpha^2-1}.
$$

Thus, with a probability of at least $1-\delta$,

\begin{align*}
\ell &\leq  \frac{\sigma^2}{\alpha-1}\left [1 + \sqrt{\frac{ 3 \alpha^2-1}{\delta}}\right].
\end{align*}

\subsection{$n/m \rightarrow \infty$}

It is easy to see that $\lim_{\alpha \rightarrow \infty} \mu_\ell  =    \sigma^2 $ and $\lim_{\alpha \rightarrow \infty} \sigma^2_\ell  =   3 \sigma^4 $. Therefore as $\alpha \rightarrow \infty$, with a probability of at least $1-\delta$, one can guarantee that

\begin{align*}
\ell &\leq  \sigma^2  \sqrt{\frac{ 3}{\delta}}.
\end{align*}

\bibliographystyle{IEEEtran}{}
\bibliography{refs}

\end{document}

%% file: preamble.tex
\usepackage[margin=1in]{geometry}
\usepackage{amsmath}
\usepackage{amssymb}
\usepackage{hyperref}
\usepackage{subcaption}
\usepackage{tikz}
\usepackage{pgfplots}

\usepackage{amsthm}
\newtheorem{theorem}{Theorem}

%% file: macros.tex
\newcommand{\reals}{\mathbb{R}}

\newcommand{\norm}[1]{\| #1 \|_2}
\newcommand{\tr}{\textrm{Tr}}
\newcommand{\xh}{\hat{x}}
\newcommand{\yh}{\hat{y}}
\newcommand{\zh}{\hat{z}}

\usepackage{empheq}
\usepackage[most]{tcolorbox}

\newtcbox{\mymath}[1][]{%
	nobeforeafter, math upper, tcbox raise base,
	enhanced, colframe=blue!30!black,
	colback=blue!30, boxrule=1pt,
	#1}

%% file: LeastSquares.bbl
\begin{thebibliography}{10}
\providecommand{\url}[1]{#1}
\csname url@samestyle\endcsname
\providecommand{\newblock}{\relax}
\providecommand{\bibinfo}[2]{#2}
\providecommand{\BIBentrySTDinterwordspacing}{\spaceskip=0pt\relax}
\providecommand{\BIBentryALTinterwordstretchfactor}{4}
\providecommand{\BIBentryALTinterwordspacing}{\spaceskip=\fontdimen2\font plus
\BIBentryALTinterwordstretchfactor\fontdimen3\font minus
  \fontdimen4\font\relax}
\providecommand{\BIBforeignlanguage}[2]{{%
\expandafter\ifx\csname l@#1\endcsname\relax
\typeout{** WARNING: IEEEtran.bst: No hyphenation pattern has been}%
\typeout{** loaded for the language `#1'. Using the pattern for}%
\typeout{** the default language instead.}%
\else
\language=\csname l@#1\endcsname
\fi
#2}}
\providecommand{\BIBdecl}{\relax}
\BIBdecl

\bibitem{breiman1983many}
L.~Breiman and D.~Freedman, ``How many variables should be entered in a
  regression equation?'' \emph{Journal of the American Statistical
  Association}, vol.~78, no. 381, pp. 131--136, 1983.

\bibitem{hastie2019surprises}
T.~Hastie, A.~Montanari, S.~Rosset, and R.~J. Tibshirani, ``Surprises in
  high-dimensional ridgeless least squares interpolation,'' \emph{arXiv
  preprint arXiv:1903.08560}, 2019.

\bibitem{gyorfi2002distribution}
L.~Gy{\"o}rfi, M.~Kohler, A.~Krzyzak, H.~Walk \emph{et~al.}, \emph{A
  distribution-free theory of nonparametric regression}.\hskip 1em plus 0.5em
  minus 0.4em\relax Springer, 2002, vol.~1.

\bibitem{catoni2004statistical}
O.~Catoni, \emph{Statistical learning theory and stochastic optimization: Ecole
  d'Et{\'e} de Probabilit{\'e}s de Saint-Flour, XXXI-2001}.\hskip 1em plus
  0.5em minus 0.4em\relax Springer Science \& Business Media, 2004, vol. 1851.

\bibitem{audibert2010linear}
J.-Y. Audibert and O.~Catoni, ``Linear regression through pac-bayesian
  truncation,'' \emph{arXiv preprint arXiv:1010.0072}, 2010.

\bibitem{hsu2012random}
D.~Hsu, S.~M. Kakade, and T.~Zhang, ``Random design analysis of ridge
  regression,'' in \emph{Conference on learning theory}.\hskip 1em plus 0.5em
  minus 0.4em\relax JMLR Workshop and Conference Proceedings, 2012.

\bibitem{mcallester1999pac}
D.~A. McAllester, ``Pac-bayesian model averaging,'' in \emph{Proceedings of the
  twelfth annual conference on Computational learning theory}, 1999, pp.
  164--170.

\bibitem{germain2016pac}
P.~Germain, F.~Bach, A.~Lacoste, and S.~Lacoste-Julien, ``Pac-bayesian theory
  meets bayesian inference,'' \emph{arXiv preprint arXiv:1605.08636}, 2016.

\bibitem{mei2019generalization}
S.~Mei and A.~Montanari, ``The generalization error of random features
  regression: Precise asymptotics and the double descent curve,''
  \emph{Communications on Pure and Applied Mathematics}, 2019.

\bibitem{gelman2013bayesian}
A.~Gelman, J.~B. Carlin, H.~S. Stern, D.~B. Dunson, A.~Vehtari, and D.~B.
  Rubin, \emph{Bayesian data analysis}.\hskip 1em plus 0.5em minus 0.4em\relax
  CRC press, 2013.

\bibitem{matrixref}
M.~Brookes, ``The matrix reference manual,''
  \emph{http://www.ee.imperial.ac.uk/hp/staff/dmb/matrix/intro.html}, 2020.

\bibitem{pielaszkiewicz2019mixtures}
J.~Pielaszkiewicz and T.~Holgersson, ``Mixtures of traces of wishart and
  inverse wishart matrices,'' \emph{Communications in Statistics-Theory and
  Methods}, pp. 1--17, 2019.

\end{thebibliography}
